\input amstex
\input amsppt.sty   
\hsize 13cm
\vsize 19.2cm
\magnification=\magstep1
\def\nmb#1#2{#2}         
\def\totoc{}             
\def\ign#1{}             

\redefine\o{\circ}
\define\X{\frak X}
\define\al{\alpha}

\define\ga{\gamma}
\define\de{\delta}

\define\la{\lambda}
\define\rh{\rho}
\define\si{\sigma}

\define\om{\omega}

\define\La{\Lambda}

\define\Om{\Omega}
\redefine\i{^{-1}}
\define\row#1#2#3{#1_{#2},\ldots,#1_{#3}}
\define\x{\times}
\redefine\L{{\Cal L}}
\define\g{{\frak g}}
\define\ome{\varOmega}           
\define\grad{\operatorname{grad}^\ome }
\define\a{\omega}                 
\define\Lip{\operatorname{Lip}}
\define\ev{\operatorname{ev}}
\def\today{\ifcase\month\or
 January\or February\or March\or April\or May\or June\or
 July\or August\or September\or October\or November\or December\fi
 \space\number\day, \number\year}
\topmatter
\title The Moment Mapping for Unitary Representations  \endtitle
\author  Peter W. Michor  \endauthor
\affil Institut f\"ur Mathematik der Universit\"at Wien, Austria \endaffil
\address{Institut f\"ur Mathematik, Universit\"at Wien,
Strudlhofgasse 4, A-1090 Wien, Austria. 
}\endaddress
\email Peter.Michor\@esi.ac.at\endemail
\abstract{For any unitary representation of an arbitrary Lie
group I construct a moment mapping from the space of smooth
vectors of the representation into the dual of the Lie algebra.
This moment mapping is equivariant and smooth. For the space of
analytic vectors the same construction is possible and leads to
a real analytic moment mapping.}\endabstract
\endtopmatter

\document
\heading Table of contents\endheading
\noindent 1. Introduction \leaders \hbox to 1em{\hss .\hss }\hfill {\eightrm 1}\par 
\noindent 2. Calculus of smooth mappings \leaders \hbox to 1em{\hss .\hss }\hfill {\eightrm 2}\par 
\noindent 3. Calculus of holomorphic mappings \leaders \hbox to 1em{\hss .\hss }\hfill {\eightrm 5}\par 
\noindent 4. Calculus of real analytic mappings \leaders \hbox to 1em{\hss .\hss }\hfill {\eightrm 6}\par 
\noindent 5. The Space of Smooth Vectors \leaders \hbox to 1em{\hss .\hss }\hfill {\eightrm 7}\par 
\noindent 6. The model for the moment mapping \leaders \hbox to 1em{\hss .\hss }\hfill {\eightrm 8}\par 
\noindent 7. Hamiltonian Mechanics on $H_{\infty }$ \leaders \hbox to 1em{\hss .\hss }\hfill {\eightrm 9}\par 
\noindent 8. The moment mapping for a unitary representation \leaders \hbox to 1em{\hss .\hss }\hfill {\eightrm 11}\par 
\noindent 9. The real analytic moment mapping \leaders \hbox to 1em{\hss .\hss }\hfill {\eightrm 13}\par 


\heading\totoc \nmb0{1}. Introduction \endheading

With the help of the cartesian closed calculus for smooth
mappings as explained in \cite{F-K} we can show, that for any
Lie group and for any unitary representation its restriction 
to the space of smooth vectors is smooth. The imaginary part of
the hermitian inner product restricts to a "weak" symplectic
structure on the vector space of smooth vectors. This gives rise
to the Poisson bracket on a suitably chosen space of smooth
functions on the space of smooth vectors. The derivative of the
representation on the space of smooth vectors is a symplectic
action of the Lie algebra, which can be lifted to a Hamiltonian
action, i.e. a Lie algebra homomorphism from the Lie algebra
into the function space with the Poisson bracket.
This in turn gives rise to the moment mapping from the space
of smooth vectors into the dual of the Lie algebra.

In \cite{K-M} the cartesian closed setting for real analytic
mappings in infinite dimensions is fully developed. With its
help it can be shown that the moment mapping
restricts to a real analytic mapping from the subspace of
analytic vectors into the dual of the Lie algebra.

For an irreducible representation which
is constructed by geometric quantization of an coadjoint orbit
(the Kirillov method), the restriction of the moment mapping to
the intersection of the unit sphere with the space of smooth
vectors takes values has as image exactly the convex hull of 
the orbit one started with, if
the construction is suitably normalized. This has been proved by 
Wildberger \cite{Wil}. I thank J. Hilgert for bringing this paper to 
my attention.

Let me add some thoughts on the r\^ole of the moment mapping in
the study of unitary representations. I think that its
restriction to the intersection of the unit sphere with the
space of smooth vectors maps to the convex hull of 
one coadjoint orbit, if the
representation is irreducible (I was unable to prove this). It
is known that not all irreducible representations come from line
bundles over coadjoint orbits (alias geometric quantization),
but there might be a higher dimensional vector bundle over this
coadjoint orbit, whose space of sections contains the space of
smooth vectors as subspace of sections which are covariantly
constant along some complex polarization.

For the convenience of the reader I have added three sections on the 
smooth, holomorphic and real analytic setting, which is used in the 
rest of the paper. These sections are of review character.

\heading\totoc\nmb0{2}. Calculus of smooth mappings \endheading

\subheading{\nmb.{2.1}} The traditional differential calculus works 
well for finite dimensional vector spaces and for Banach spaces. For 
more general locally convex spaces a whole flock of different 
theories were developed, each of them rather complicated and none 
really convincing. The main difficulty is that the composition of 
linear mappings stops to be jointly continuous at the level of Banach 
spaces, for any compatible topology. This was the original motivation 
for the development of a whole new field within general topology, 
convergence spaces.

Then in 1982, Alfred Fr\"olicher and Andreas Kriegl presented 
independently the solution to the question for the right differential 
calculus in infinite dimensions. They joined forces in the further 
development of the theory and the (up to now) final outcome is the 
book \cite{F-K}.

In this section I will sketch the basic definitions and the most 
important results of the Fr\"olicher-Kriegl calculus.

\subheading{\nmb.{2.2}. The $c^\infty$-topology} Let $E$ be a 
locally convex vector space. A curve $c:\Bbb R\to E$ is called 
{\it smooth} or $C^\infty$ if all derivatives exist and are 
continuous - this is a concept without problems. Let 
$C^\infty(\Bbb R,E)$ be the space of smooth functions. It can be 
shown that $C^\infty(\Bbb R,E)$ does not depend on the locally convex 
topology of $E$, only on its associated bornology (system of bounded 
sets).

The final topologies with respect to the following sets of mappings 
into E coincide:
\roster
\item $C^\infty(\Bbb R,E)$.
\item Lipschitz curves (so that $\{\frac{c(t)-c(s)}{t-s}:t\neq s\}$ 
     is bounded in $E$). 
\item $\{E_B\to E: B\text{ bounded absolutely convex in }E\}$, where 
     $E_B$ is the linear span of $B$ equipped with the Minkowski 
     functional $p_B(x):= \inf\{\la>0:x\in\la B\}$.
\item Mackey-convergent sequences $x_n\to x$ (there exists a sequence 
     $0<\la_n\nearrow\infty$ with $\la_n(x_n-x)$ bounded).
\endroster
This topology is called the $c^\infty$-topology on $E$ and we write 
$c^\infty E$ for the resulting topological space. In general 
(on the space $\Cal D$ of test functions for example) it is finer 
than the given locally convex topology, it is not a vector space 
topology, since scalar multiplication is no longer jointly 
continuous. The finest among all locally convex topologies on $E$ 
which are coarser than $c^\infty E$ is the bornologification of the 
given locally convex topology. If $E$ is a Fr\'echet space, then 
$c^\infty E = E$. 

\subheading{\nmb.{2.3}. Convenient vector spaces} Let $E$ be a 
locally convex vector space. $E$ is said to be a {\it convenient 
vector space} if one of the following equivalent (completeness) 
conditions is satisfied:
\roster
\item Any Mackey-Cauchy-sequence (so that $(x_n-x_m)$ is Mackey 
     convergent to 0) converges. This is also called 
     $c^\infty$-complete.
\item If $B$ is bounded closed absolutely convex, then $E_B$ is a 
     Banach space.
\item Any Lipschitz curve in $E$ is locally Riemann integrable.
\item For any $c_1\in C^\infty(\Bbb R,E)$ there is 
     $c_2\in C^\infty(\Bbb R,E)$ with $c_1'=c_2$ (existence of 
     antiderivative).
\endroster

\proclaim{\nmb.{2.4}. Lemma} Let $E$ be a locally convex space.
Then the following properties are equivalent:
\roster
\item $E$ is $c^\infty$-complete.
\item If $f:\Bbb R^k\to E$ is scalarwise $\Lip^k$, then $f$ is 
     $\Lip^k$, for $k>1$.
\item If $f:\Bbb R\to E$ is scalarwise $C^\infty$ then $f$ is 
     differentiable at 0.
\item If $f:\Bbb R\to E$ is scalarwise $C^\infty$ then $f$ is 
     $C^\infty$.
\endroster
\endproclaim
Here a mapping $f:\Bbb R^k\to E$ is called $\Lip^k$ if all partial 
derivatives up to order $k$ exist and are Lipschitz, locally on 
$\Bbb R^n$. $f$ scalarwise $C^\infty$ means that $\la\o f$ is $C^\infty$  
for all continuous linear functionals on $E$.

This lemma says that a convenient vector space one can recognize 
smooth curves by investigating compositions with continuous linear 
functionals.

\subheading{\nmb.{2.5}. Smooth mappings} Let $E$ and $F$ be locally 
convex vector spaces. A mapping $f:E\to F$ is called {\it smooth} or 
$C^\infty$, if $f\o c\in C^\infty(\Bbb R,F)$ for all 
$c\in C^\infty(\Bbb R,E)$; so 
$f_*: C^\infty(\Bbb R,E)\to C^\infty(\Bbb R,F)$ makes sense.
Let $C^\infty(E,F)$ denote the space of all smooth mapping from $E$ 
to $F$.

For $E$ and $F$ finite dimensional this gives the usual notion of 
smooth mappings: this has been first proved in \cite{Bo}.
Constant mappings are smooth. Multilinear mappings are smooth if and 
only if they are boun\-ded. Therefore we denote by $L(E,F)$ the space 
of all bounded linear mappings from $E$ to $F$.

\subheading{\nmb.{2.6}. Structure on $C^\infty(E,F)$} We equip the 
space $C^\infty(\Bbb R,E)$ with the bornologification of the topology 
of uniform convergence on compact sets, in all derivatives 
separately. Then we equip the space $C^\infty(E,F)$ with the 
bornologification of the initial topology with respect to all 
mappings $c^*:C^\infty(E,F)\to C^\infty(\Bbb R,F)$, $c^*(f):=f\o c$, 
for all $c\in C^\infty(\Bbb R,E)$.

\proclaim{\nmb.{2.7}. Lemma } For locally convex spaces $E$ and $F$ 
we have:
\roster
\item If $F$ is convenient, then also $C^\infty(E,F)$ is convenient, 
     for any $E$. The space $L(E,F)$ is a closed linear subspace of 
     $C^\infty(E,F)$, so it also  convenient.
\item If $E$ is convenient, then a curve $c:\Bbb R\to L(E,F)$ is 
     smooth if and only if $t\mapsto c(t)(x)$ is a smooth curve in $F$ 
     for all $x\in E$.
\endroster
\endproclaim

\proclaim{\nmb.{2.8}. Theorem} The category of convenient vector 
spaces and smooth mappings is cartesian closed. So we have a natural 
bijection 
$$C^\infty(E\x F,G)\cong C^\infty(E,C^\infty(F,G)),$$
which is even a diffeomorphism.
\endproclaim

Of coarse this statement is also true for $c^\infty$-open subsets of 
convenient vector spaces. 

\proclaim{\nmb.{2.9}. Corollary } Let all spaces be convenient vector 
spaces. Then the following canonical mappings are smooth.
$$\align
&\operatorname{ev}: C^\infty(E,F)\x E\to F,\quad 
     \operatorname{ev}(f,x) = f(x)\\
&\operatorname{ins}: E\to C^\infty(F,E\x F),\quad
     \operatorname{ins}(x)(y) = (x,y)\\
&(\quad)^\wedge :C^\infty(E,C^\infty(F,G))\to C^\infty(E\x F,G)\\
&(\quad)\spcheck :C^\infty(E\x F,G)\to C^\infty(E,C^\infty(F,G))\\
&\operatorname{comp}:C^\infty(F,G)\x C^\infty(E,F)\to C^\infty(E,G)\\
&C^\infty(\quad,\quad):C^\infty(F,F')\x C^\infty(E',E)\to 
     C^\infty(C^\infty(E,F),C^\infty(E',F'))\\
&\qquad (f,g)\mapsto(h\mapsto f\o h\o g)\\
&\prod:\prod C^\infty(E_i,F_i)\to C^\infty(\prod E_i,\prod F_i)
\endalign$$
\endproclaim

\proclaim{\nmb.{2.10}. Theorem} Let $E$ and $F$ be convenient vector 
spaces. Then the differential operator 
$$\gather d: C^\infty(E,F)\to C^\infty(E,L(E,F)), \\
df(x)v:=\lim_{t\to0}\frac{f(x+tv)-f(x)}t,
\endgather$$
exists and is linear and bounded (smooth). Also the chain rule holds: 
$$d(f\o g)(x)v = df(g(x))dg(x)v.$$
\endproclaim

\subheading{\nmb.{2.11}. Remarks } Note that the conclusion of 
theorem \nmb!{2.8} is the starting point of the classical calculus of 
variations, where a smooth curve in a space of functions was assumed 
to be just a smooth function in one variable more.

If one wants theorem \nmb!{2.8} to be true and assumes some other obvious 
properties, then the calculus of smooth functions is already uniquely 
determined.

There are, however, smooth mappings which are not continuous. This is 
unavoidable and not so horrible as it might appear at first sight. 
For example the evaluation $E\x E'\to\Bbb R$ is jointly continuous if 
and only if $E$ is normable, but it is always smooth. Clearly smooth 
mappings are continuous for the $c^\infty$-topology.

For Fr\'echet spaces smoothness in the sense described here coincides 
with the notion $C^\infty_c$ of \cite{Ke}. This is the 
differential calculus used by \cite{Mic}, \cite{Mil}, 
and \cite{P-S}.

\heading\totoc\nmb0{3}. Calculus of holomorphic mappings \endheading

\subheading{\nmb.{3.1}} Along the lines of thought of the 
Fr\"olicher-Kriegl calculus of smooth mappings, in \cite{K-N} 
the cartesian closed setting for holomorphic mappings was 
developed. 
We will now sketch the basics 
and the main results. It can be shown that again convenient vector 
spaces are the right ones to consider. Here we will start with them 
for the sake of shortness.

\subheading{\nmb.{3.2}} Let $E$ be a complex locally convex vector 
space whose underlying real space is convenient -- this will be 
called convenient in the sequel. Let $\Bbb D\subset \Bbb C$ be the 
open unit disk and let us denote by $\Cal H(\Bbb D,E)$ the space of 
all mappings $c:\Bbb D\to E$ such that $\la\o c:\Bbb D\to \Bbb C$ is 
holomorphic for each continuous complex-linear functional $\la$ on 
$E$. Its elements will be called the holomorphic curves.

If $E$ and $F$ are convenient complex vector spaces (or 
$c^\infty$-open sets therein), a mapping 
$f:E\to F$ is called {\it holomorphic} if $f\o c$ is a holomorphic 
curve in $F$ for each holomorphic curve $c$ in $E$. Obviously $f$ is 
holomorphic if and only if $\la\o f:E\to \Bbb C$ is holomorphic for 
each complex linear continuous functional $\la$ on $F$. Let 
$\Cal H(E,F)$ denote the space of all holomorphic mappings from $E$ to 
$F$. 

\proclaim{\nmb.{3.3}. Theorem (Hartog's theorem)} Let $E_k$ for 
$k=1,2$ and $F$ be complex convenient vector spaces and let 
$U_k\subset E_k$ be $c^\infty$-open. A mapping $f:U_1\x U_2\to F$ is 
holomorphic if and only if it is separately holomorphic (i\. e\. 
$f(\quad,y)$ and $f(x,\quad)$ are holomorphic for all $x\in U_1$ and 
$y\in U_2$).
\endproclaim

This implies also that in finite dimensions we have recovered the 
usual definition.

\proclaim{\nmb.{3.4} Lemma} If $f:E\supset U\to F$ is holomorphic 
then $df:U\x E\to F$ exists, is holomorphic and $\Bbb C$-linear in 
the second variable. 

A multilinear mapping is holomorphic if and only if it is bounded.
\endproclaim

\proclaim{\nmb.{3.5} Lemma} If $E$ and $F$ are Banach spaces and $U$ 
is open in $E$, then for a mapping $f:U\to F$ the following 
conditions are equivalent:
\roster
\item $f$ is holomorphic.
\item $f$ is locally a convergent series of homogeneous continuous 
     polynomials.
\item $f$ is $\Bbb C$-differentiable in the sense of Fr\'echet.
\endroster 
\endproclaim

\proclaim{\nmb.{3.6} Lemma} Let $E$ and $F$ be convenient vector 
spaces. A mapping $f:E\to F$ is holomorphic if and only if it is 
smooth and its derivative is everywhere $\Bbb C$-linear.
\endproclaim

An immediate consequence of this result is that $\Cal H(E,F)$ is a 
closed linear subspace of $C^\infty(E_{\Bbb R},F_{\Bbb R})$ and so it 
is a convenient vector space if $F$ is one, by \nmb!{2.7}. The chain 
rule follows from \nmb!{2.10}. The following theorem is an easy 
consequence of \nmb!{2.8}.

\proclaim{\nmb.{3.7} Theorem} The category of convenient complex 
vector spaces and holomorphic mappings between them is cartesian 
closed, i\. e\.
$$\Cal H(E\x F,G) \cong \Cal H(E,\Cal H(F,G)).$$
\endproclaim

An immediate consequence of this is again that all canonical 
structural mappings as in \nmb!{2.9} are holomorphic.

\heading\totoc\nmb0{4}. Calculus of real analytic mappings \endheading 

\subheading{\nmb.{4.1}} In this section we sketch the cartesian closed 
setting to real analytic mappings in infinite dimension following the 
lines of the Fr\"olicher-Kriegl calculus, as it is presented in 
\cite{K-M}. Surprisingly enough one has to deviate 
from the most obvious notion of real analytic curves in order to get 
a meaningful theory, but again convenient vector spaces turn out to 
be the right kind of spaces.

\subheading{\nmb.{4.2}. Real analytic curves} Let $E$ be a real 
convenient vector space with dual $E'$. A curve $c:\Bbb R\to E$ is 
called {\it real analytic} if $\la\o c:\Bbb R\to \Bbb R$ is real 
analytic for each $\la\in E'$.
It turns out that the set of these curves depends only on the 
bornology of $E$.

In contrast a curve is called {\it topologically real analytic} if it 
is locally given by power series which converge in the topology of 
$E$. They can be extended to germs of holomorphic curves along $\Bbb R$ 
in the complexification $E_{\Bbb C}$ of $E$. If the dual $E'$ of $E$ 
admits a Baire topology which is compatible with the duality, then 
each real analytic curve in $E$ is in fact topologically real analytic 
for the bornological topology on $E$.

\subheading{\nmb.{4.3}. Real analytic mappings} Let $E$ and $F$ be 
convenient vector spaces. Let $U$ be a $c^\infty$-open set in $E$. A 
mapping $f:U\to F$ is called {\it real analytic} if and only if it is 
smooth (maps smooth curves to smooth curves) and maps real analytic 
curves to real analytic curves. 

Let $C^\omega(U,F)$ denote the space of all real analytic mappings. 
We equip the space
$C^\om(U,\Bbb R)$ of all real analytic functions 
with the initial topology
with respect to the families of mappings
$$\gather C^\om(U,\Bbb R) @>{c^*}>> C^\om(\Bbb R,\Bbb R),\text{
for all }c\in C^\om(\Bbb R,U)\\ 
C^\om(U,\Bbb R) @>{c^*}>> C^\infty(\Bbb R,\Bbb R),\text{
for all }c\in C^\infty(\Bbb R,U),
\endgather$$
where $C^\infty(\Bbb R,\Bbb R)$ carries the topology of compact 
convergence in each derivative separately as in section \nmb!{2}, and 
where $C^\omega(\Bbb R,\Bbb R)$ is equipped with the final locally 
convex topology 
with respect to the embeddings (restriction mappings) of all spaces 
of holomorphic mappings from a neighborhood $V$ of $\Bbb R$ in 
$\Bbb C$ mapping $\Bbb R$ to $\Bbb R$, and each of these spaces 
carries the topology of compact convergence.

Furthermore we equip the space
$C^\om(U,F)$ with the initial topology with respect to
the family of mappings 
$$ C^\om(U,F) @>{\la_*}>> C^\om(U,\Bbb R),\text{ for all }\la\in
F'.$$
It turns out that this is again a convenient space.

\proclaim{\nmb.{4.4}. Theorem} In the setting of \nmb!{4.3} a mapping 
$f:U\to F$ is real analytic if and only if it is smooth and is real 
analytic along each affine line in $E$.
\endproclaim

\proclaim{\nmb.{4.5}. Lemma}The space $L(E,F)$ of all bounded linear 
mappings is a closed linear subspace of $C^\omega(E,F)$. A mapping 
$f:U\to L(E,F)$ is real analytic if and only if $\ev_x\o f:U\to F$ is 
real analytic for each point $x\in E$.
\endproclaim

\proclaim{\nmb.{4.6}. Theorem} The category of convenient spaces and 
real analytic mappings is cartesian closed. So the equation
$$C^\omega(U,C^\omega(V,F))\cong C^\omega(U\x V,F)$$
is valid for all $c^\infty$-open sets $U$ in $E$ and $V$ in $F$, 
where $E$, $F$, and $G$ are convenient vector spaces.
\endproclaim

This implies again that all structure mappings as in \nmb!{2.9} are 
real analytic. Furthermore the differential operator 
$$d:C^\omega(U,F)\to C^\omega(U,L(E,F))$$ 
exists, is unique and real 
analytic. Multilinear mappings are real analytic if and only if they 
are bounded. Powerful real analytic uniform boundedness principles 
are available.

\heading\totoc \nmb0{5}. The Space of Smooth Vectors \endheading

\subheading{\nmb.{5.1}} Let $G$ be any (finite dimensional second
countable) real Lie group, and let $\rh: G\to  U(\bold H)$ be a
unitary representation on a Hilbert space $\bold H$. Then the
associated mapping $\hat \rh: G\x \bold H \to  \bold H$ is in general
{\it not} jointly continuous, it is only separately continuous,
so that $g\mapsto \rh(g)x$, $G\to  \bold H$, is continuous for any
$x\in \bold H$.

\subheading{Definition} A vector $x \in \bold H$ is called {\it
smooth} (or {\it real analytic}) if the mapping $g\mapsto
\rh(g)x$, $G\to \bold H$ is smooth (or real analytic).
Let us denote by $\bold H_\infty$ the linear subspace of all
smooth vectors in $\bold H$. Then we have an embedding 
$j:\bold H_\infty \to  C^\infty(G,\bold H)$, given by 
$x\mapsto(g\mapsto\rh(g)x)$. We equip $C^\infty(G,\bold H)$ with
the compact $C^\infty$-topology (of uniform convergence on
compact subsets of $G$, in all derivatives separately). Then
it is easily seen (and proved in \cite {Wa, p 253}) that 
$\bold H_\infty$ is a closed linear subspace. So with the
induced topology $\bold H_\infty$ becomes a Fr\`echet space.
Clearly $\bold H_\infty$ is also an invariant subspace, so we
have a representation $\rh: G\to  L(\bold H_\infty,\bold H_\infty)$.
For more detailed information on $\bold H_\infty$ see
\cite{Wa, chapt. 4.4.} or \cite{Kn, chapt. III.}.

\proclaim{\nmb.{5.2}. Theorem} The mapping $\hat\rh: G\x\bold H_\infty \to 
\bold H_\infty$ is smooth in the sense of Fr\"olicher-Kriegl.
\endproclaim
\demo{Proof} By cartesian closedness \nmb!{2.8}
it suffices to
show that the canonically associated mapping 
$$\hat\rh\spcheck :G \to  C^\infty(\bold H_\infty,\bold H_\infty)$$
is smooth; but it takes values in the closed subspace 
$L(\bold H_\infty,\bold H_\infty)$ 
of all bounded linear operators. So by  it suffices to show that the mapping 
$\rh:G \to  L(\bold H_\infty,\bold H_\infty)$
is smooth. But for that, since $\bold H_\infty$ is a Fr\`echet
space and thus convenient,
by \nmb!{2.7}\therosteritem2 it suffices to show that 
$$G @>{\rho}>>  L(\bold H_\infty,\bold H_\infty) @>{ev_x}>> \bold H_\infty$$ 
is smooth for each $x \in \bold
H_\infty$. This requirement means that $g\mapsto \rh(g)x$,
$G\to \bold H_\infty$, is smooth.  For this it suffices to show
that $$\align &G\to \bold H_\infty @>j>> C^\infty(G,\bold H),\\
&g\mapsto \rh(g)x \mapsto(h\mapsto \rh(h)(g)x),\endalign$$ is
smooth.  But again by cartesian closedness it suffices to show
that the associated mapping $$\align &G\x G\to  \bold H,\\
&(g,h)\mapsto \rh(h)(g)x = \rh(hg)x,\endalign$$ is smooth. And
this is the case since $x$ is a smooth vector.
\qed\enddemo

\heading\totoc \nmb0{6}. The model for the moment mapping \endheading

\subheading{\nmb.{6.1}} We now consider $\bold H_\infty$ as a "weak"
symplectic Fr\`echet manifold, equipped with the symplectic
structure $\ome$, the restriction of the imaginary part of the
Hermitian inner product $\langle \quad,\quad\rangle$ on $\bold H$. Then
$\ome\in \Om^2(\bold H_\infty)$ is a closed 2-form which is non
degenerate in the sense that 
$$\check \ome: T\bold H_\infty = \bold H_\infty\x\bold H_\infty \to 
T^*\bold H_\infty = \bold H_\infty\x{\bold H_\infty}'$$
is injective (but not surjective), where 
$\bold H_\infty{}' = L(\bold H_\infty,\Bbb R)$
denotes the real topological dual space. This is the meaning of "weak"
above.

\subheading{\nmb.{6.2}. Review} For a finite dimensional symplectic
manifold $(M,\ome)$ we have the following exact sequence of
Lie algebras:
$$0\to  H^0(M)\to C^\infty(M)  @>{\grad}>> \X_\ome(M) @>\ga>> H^1(M)\to 0$$
Here $H^*(M)$ is the real De Rham cohomology of $M$, the space
$C^\infty(M)$ is equipped with the Poisson bracket $\{\quad,\quad\}$,
$\X_\ome(M)$ consists of all vector fields $\xi$ with $\L_\xi\ome=0$
(the locally Hamiltonian vector fields), which is a Lie algebra
for the Lie bracket. Also $\grad f$ is the Hamiltonian vector field
for $f\in C^\infty(M)$ given by $i(\grad f)\ome = df$, and 
$\ga(\xi) = [i_\xi\ome]$. The spaces $H^0(M)$ and $H^1(M)$ are
equipped with the zero bracket.

Given a symplectic left action $\ell:G\x M\to  M$ of a connected
Lie group $G$ on $M$, the first partial derivative of $\ell$
gives a mapping $\ell':\g \to  \X_\ome(M)$ which sends each element
$X$ of the Lie algebra $\g$ of $G$ to the fundamental vector
field. This is a Lie algebra homomorphism. 
$$\CD
H^0(M) @>i>> C^\infty(M)  @>{\grad}>> \X_\ome(M) @>\ga>> H^1(M) \\
@.           @A\si AA            @AA{\ell'}A      @.  \\
       @.    \g           @=    \g        @.            
\endCD$$
A linear lift $\si:\g\to  C^\infty(M)$ of $\ell'$  with 
$\grad\o\si=\ell'$ exists if
and only if $\ga\o\ell'=0$ in $H^1(M)$. This lift $\si$ may be
changed to a Lie algebra homomorphism if and only if the
$2$-cocycle $\bar\si:\g\x\g\to  H^0(M)$, given by 
$(i\o\bar\si)(X,Y) = \{\si(X),\si(Y)\} - \si([X,Y])$, vanishes
in $H^2(\g,H^0(M))$, for if $\bar\si = \de\al$ then $\si-i\o\al$
is a Lie algebra homomorphism.

If $\si:\g\to C^\infty(M)$ is a Lie algebra homomorphism, we may
associate the {\it moment mapping} 
$\mu:M\to \g'=L(\g,\Bbb R)$ to it, which is 
given by
$\mu(x)(X) = \si(X)(x)$ for $x\in M$ and $X\in \g$.
It is $G$-equivariant for a suitably chosen (in general affine)
action of $G$ on $\g'$.
See \cite{We} or \cite{L-M} for all this.

\heading\totoc \nmb0{7}. Hamiltonian Mechanics on $\bold H_{\infty}$ 
\endheading

\subheading{\nmb.{7.1}} We now want to carry over to the setting of \nmb!{5.1}
and \nmb!{5.2} the procedure of \nmb!{6.2}. The first thing to note is that
the hamiltonian mapping 
$\grad:C^\infty(\bold H_\infty) \to  \X_\ome(\bold H_\infty)$
does not make sense in general, since 
$\check\ome:\bold H_\infty \to  {\bold H_\infty}'$ is not invertible:
$\grad f = \check\ome\i df$ is defined only for those $f\in
C^\infty(\bold H_\infty)$ with $df(x)$ in the image of $\check\ome$ for all
$x\in \bold H_\infty$. A similar difficulty arises for the
definition of the Poisson bracket on $C^\infty(\bold H_\infty)$.

Let $\langle x,y\rangle = Re\langle x,y\rangle +
\sqrt{-1}\ome(x,y)$ be the decomposition of the hermitian inner
product into real and imaginary parts. Then $Re\langle
x,y\rangle = \ome(\sqrt{-1}x,y)$, thus the real linear subspaces
$\check\ome(\bold H_\infty) = \ome(\bold H_\infty,\quad)$ and 
$Re\langle \bold H_\infty,\quad\rangle$ of $\bold H_\infty{}' =
L(\bold H_\infty,\Bbb R)$ coincide.

\subheading{\nmb.{7.2} Definition} Let $\bold H_\infty^*$ denote
the real linear subspace 
$$\bold H_\infty^*=\ome(\bold H_\infty,\quad) = 
Re\langle \bold H_\infty,\quad\rangle$$ 
of $\bold H_\infty{}' = 
L(\bold H_\infty,\Bbb R)$, and let us call it the {\it smooth dual}
of $\bold H_\infty$ in view of the embedding of test functions
into distributions. We have two canonical isomorphisms 
$\bold H_\infty^* \cong \bold H_\infty$ induced by $\ome$ and 
$Re\langle\quad,\quad\rangle$, respectively. Both induce the
same Fr\'echet topology on $\bold H_\infty^*$, which we fix from
now on.

\subheading{\nmb.{7.3} Definition} Let 
$C^\infty_*(\bold H_\infty,\Bbb R) \subset
C^\infty(\bold H_\infty,\Bbb R)$ 
denote the linear subspace consisting of all smooth functions
$f:\bold H_\infty \to  \Bbb R$ such that each iterated derivative 
$d^kf(x)\in L^k_{\text{sym}}(\bold H_\infty,\Bbb R)$ has the
property that 
$$d^kf(x)(\quad,\row y2k) \in {\bold H_\infty}^*$$
is actually in the smooth dual $\bold H_\infty^* \subset
{\bold H_\infty}'$ for all $x,\row y2k \in \bold H_\infty$, and
that the mapping 
$$\gather \prod^k\bold H_\infty \to  \bold H_\infty \\
(x,\row y2k)\mapsto \check\ome\i(df(x)(\quad,\row y2k))
\endgather$$
is smooth. Note that we could also have used
$Re\langle\quad,\quad\rangle$ instead of $\ome$. By the symmetry
of higher derivatives this is then true for all entries of
$d^kf(x)$, for all $x$.  

\proclaim{\nmb.{7.4} Lemma} For $f \in C^\infty(\bold H_\infty,\Bbb R)$ the
following assertions are equivalent:
\roster
\item $df:\bold H_\infty \to  \bold H_\infty{}'$ factors to a
        smooth mapping $\bold H_\infty \to  \bold H_\infty^*$.
\item $f$ has a smooth $\ome$-gradient 
        $\grad f\in \X(\bold H_\infty) = C^\infty(\bold
        H_\infty,\bold H_\infty)$ such that 
        $df(x)y = \ome(\grad f(x),y)$.
\item $f \in C^\infty_*(\bold H_\infty,\Bbb R)$.
\endroster
\endproclaim
\demo{Proof} Clearly \therosteritem3 $\Longrightarrow$ \therosteritem2
$\Longleftrightarrow$ \therosteritem1. We have to show that
\therosteritem2 $\Longrightarrow$ \therosteritem3.

Suppose that $f:\bold H_\infty \to  \Bbb R$ is smooth and
$df(x)y = \ome(\grad f(x),y)$. Then 
$$\align
d^kf(x)(\row y1k) &= d^kf(x)(\row y2k,y_1)  \\
&= (d^{k-1}(df)(x)(\row y2k)(y_1)  \\
&= \ome(d^{k-1}(\grad f)(x)(\row y2k),y_1).\qed
\endalign$$
\enddemo

\proclaim{\nmb.{7.5}. Theorem} 
The mapping 
$\grad :C^\infty_*(\bold H_\infty,\Bbb R) \to  \X_\ome(\bold H_\infty)$,
given by $\grad f := \check\ome\i\o df$, is well defined; also the
Poisson bracket 
$$\gather
\{\quad,\quad\}: C^\infty_*(\bold H_\infty,\Bbb R) \x
C^\infty_*(\bold H_\infty,\Bbb R) \to  
C^\infty_*(\bold H_\infty,\Bbb R),\\
\{f,g\}:= i(\grad f)i(\grad g)\ome = \ome(\grad g,\grad f) =\\
= (\grad f)(g) = dg(\grad f) 
\endgather$$ 
is well defined and gives a Lie algebra structure to
the space $C^\infty_*(\bold H_\infty,\Bbb R)$.

We also have the following long exact sequence of Lie algebras
and Lie algebra homomorphisms:
$$0\to  H^0(\bold H_\infty)\to C^\infty_*(\bold H_\infty,\Bbb R)  
@>{\grad}>> \X_\ome(\bold H_\infty) @>\ga>> H^1(\bold H_\infty) = 0$$ 
\endproclaim
\demo{Proof} It is clear from lemma \nmb!{7.4}, that the
hamiltonian mapping is defined, and thus also the Poisson
bracket is defined as a mapping
$\{\quad,\quad\}: C^\infty_*(\bold H_\infty,\Bbb R) \x
C^\infty_*(\bold H_\infty,\Bbb R) \to  
C^\infty(\bold H_\infty,\Bbb R) $,
and it only remains to check that it has values in the subspace
$C^\infty_*(\bold H_\infty,\Bbb R)$.  

So let $f$, $g\in C^\infty_*(\bold H_\infty)$, then 
$\{f,g\}(x) = dg(x)(\grad f(x))$ and by the symmetry of
$dg(x)$ we have
$$\align &d(\{f,g\})(x)y 
= d^2g(x)(y,\grad f(x)) + dg(x)(d(\grad f)(x)y)\\
&\quad = \ome\Bigl(d(\grad g)(x)(\grad f(x)),y\Bigr)\\ 
&\qquad\qquad\qquad    + \ome\Bigl(\grad g(x),d(\grad f)(x)y\Bigr) \\
&\quad = \ome\Bigl(d(\grad g)(x)(\grad f(x)) 
    - d(\grad f)(x)(\grad g(x)),y\Bigr), 
\endalign$$
since $\grad f\in\X_\ome(\bold H_\infty)$ and for any
$X\in\X_\ome(\bold H_\infty)$ the condition $\L_X\ome=0$ implies
$\ome(dX(x)y_1,y_2) = -\ome(y_1,dX(x)y_2)$. So \therosteritem2
of lemma \nmb!{7.4} is satisfied and thus 
$\{f,g\}\in C^\infty_*(\bold H_\infty)$.

For the rest any coordinate free finite dimensional proof works.
\qed \enddemo

\heading\totoc \nmb0{8}. The moment mapping for a unitary
representation \endheading

\subheading{\nmb.{8.1}} We consider now again as in
\nmb!{5.1} a unitary representation $\rh: G\to  U(\bold H)$. By theorem
\nmb!{5.2} the associated mapping 
$\hat \rh: G\x \bold H_\infty \to  \bold H_\infty$ 
is smooth, so we have the infinitesimal mapping 
$\rh':\g \to  \X(\bold H_\infty)$, given by 
$\rh'(X)(x) = T_e(\hat\rho(\quad,x))X$ for $X\in \g$ and $x\in
\bold H_\infty$. Since $\rho$ is a unitary
representation, the mapping $\rho'$ has values in the Lie
subalgebra of all linear hamiltonian vector fields $\xi \in
\X(\bold H_\infty)$ which respect the symplectic form $\ome$,
i.e. $\xi:\bold H_\infty \to  \bold H_\infty$ is linear and
$\L_\xi\ome = 0$.  

Now let us consider the mapping 
$\check\ome\o \rho'(X): \bold H_\infty \to  T(\bold H_\infty) \to 
T^*(\bold H_\infty)$. 
We have $d(\check\ome\o\rho'(X)) = d(i_{\rho'(X)}\ome) =
\L_{\rho'(X)}\ome = 0$, so the linear 1-form
$\check\ome\o\rho'(X)$ is closed, and since $H^1(\bold H_\infty)
= 0$, it is exact. So there is a function 
$\si(X)\in C^\infty(\bold H_\infty,\Bbb R)$ with $d\si(X) =
\check\ome\o\rho'(X)$, 
and $\si(X)$ is uniquely determined up to addition of a constant.
If we require $\si(X)(0) = 0$, then $\si(X)$ is uniquely
determined and is a quadratic function. In fact we have 
$\si(X)(x) = \int_{c_x}\check\ome\o\rho'(X)$, where $c_x(t) =
tx$. Thus
$$\align 
\si(X)(x) 
&= {\tsize\int_0^1} \ome(\rho'(X)(tx),\tfrac d{dt}tx)dt = \\
&= \ome(\rho'(X)(x),x){\tsize\int_0^1tdt} \\
&= \tfrac 12 \ome(\rho'(X)(x),x).
\endalign$$

\proclaim{\nmb.{8.2}. Lemma} 
The mapping 
$$\si:\g \to  C^\infty_*(\bold H_\infty,\Bbb R),\qquad 
\si(X)(x) = \tfrac12\ome(\rho'(X)(x),x)$$
for $X\in \g$ and $x\in
\bold H_\infty$, is a Lie algebra homomorphism
and $\grad\o\si = \rho'$.

For $g\in G$ we have $\rho(g)^*\si(X)=\si(X)\o \rho(g) =
\si(Ad(g\i)X)$, so $\si$ is $G$-equivariant.
\endproclaim
\demo{Proof} First we have to check that 
$\si(X)\in C^\infty_*(\bold H_\infty,\Bbb R)$. Since 
$\rho'(X):\bold H_\infty\to \bold H_\infty$ is smooth  and linear,
i.e. bounded linear, this follows from the formula for $\si(X)$.
Furthermore 
$$\align 
\grad(\si(X))(x) &= \check\ome\i(d\si(X)(x)) = \\
&= \tfrac12\check\ome\i\left(\ome(\rho'(X)(\quad),x)
+ \ome(\rho'(X)(x),\quad)\right) = \\
&=\check\ome\i\left(\ome(\rho'(X)(x),\quad)\right) = \rho'(X)(x),
\endalign$$
since $\ome(\rho'(X)(x),y) = \ome(\rho'(X)(y),x)$.

Clearly $\si([X,Y]) - \{\si(X),\si(Y)\}$ is a constant function
by \nmb!{7.5}; since it also vanishes at $0 \in \bold H_\infty$, the
mapping $\si: \g \to  C^\infty_*(\bold H_\infty)$ is a Lie
algebra homomorphism. 

For the last assertion we have
$$\align \si(X)(\rho(g)x) &= \tfrac12\ome(\rho'(X)(\rho(g)x),\rho(g)x)\\
&= \tfrac12(\rho(g)^*\ome)(\rho(g\i)\rho'(X)(\rho(g)x),x)\\
&= \tfrac12\ome(\rho'(Ad(g\i)X)x,x) = \si(Ad(g\i)X)(x).\qed
\endalign$$
\enddemo

\subheading{\nmb.{8.3}. The moment mapping} For a unitary
representation $\rho: G \to  U(\bold H)$ we can now define the
{\it moment mapping} 
$$\gather 
\mu: \bold H_\infty \to  \g' = L(\g,\Bbb R), \\
\mu(x)(X) := \si(X)(x) = \tfrac12 \ome(\rho'(X)x,x),
\endgather$$
for $x \in \bold H_\infty$ and $X \in \g$.

\proclaim{\nmb.{8.4} Theorem} The moment mapping 
$\mu:\bold H_\infty \to  \g'$ has the following properties:
\roster
\item $(d\mu(x)y)(X)=\ome(\rho'(X)x,y)$ for 
    $x,y\in \bold H_\infty$ and $X\in\g$, so 
    $\mu\in C^\infty_*(\bold H_\infty,\g')$.
\item For $x\in\bold H_\infty$ the image of 
    $d\mu(x):\bold H_\infty\to\g'$ is the
    annihilator $\g_x^\ome$ of the Lie algebra 
    $\g_x = \{X\in\g:\rho'(X)(x)=0\}$  of the isotropy group
    $G_x = \{g\in G:\rho(g)x=x\}$ in $\g'$.
\item For $x\in\bold H_\infty$ the kernel of $d\mu(x)$ is
    $$(T_x(\rho(G)x))^\ome = 
    \{y\in\bold H_\infty:\ome(y,T_x(\rho(G)x))=0\},$$
    the $\ome$-annihilator
    of the tangent space at $x$ of the $G$-orbit through $x$.
\item The moment mapping is equivariant:
    $Ad'(g)\o \mu = \mu \o \rho(g)$ for all $g\in G$,
    where $Ad'(g)=Ad(g\i)':\g'\to\g'$ is the coadjoint action.
\item The pullback operator 
    $\mu^*: C^\infty(\g,\Bbb R) \to C^\infty(\bold H_\infty,\Bbb R)$ 
    actually has values in the subspace 
    $C^\infty_*(\bold H_\infty,\Bbb R)$.
    It also is a Lie algebra homomorphism for the Poisson
    brackets involved.
\endroster
\endproclaim
\demo{Proof} \therosteritem1. Differentiating the defining
equation we get
$$(d\mu(x)y)(X) = \tfrac12\ome(\rho'(X)y,x) + \tfrac12\ome(\rho'(X)x,y) 
=\ome(\rho'(X)x,y).\tag a$$
 From lemma \nmb!{7.4} we see that $\mu\in C^\infty_*(\bold H_\infty,\g')$.

\therosteritem2 and \therosteritem3 are immediate consequences
of this formula.

\therosteritem4. We have
$$\align \mu(\rho(g)x)(X) &= \si(X)(\rho(g)x) 
    = \si(Ad(g\i)X)(x)\text{ by lemma \nmb!{8.2}}\\
&= \mu(x)(Ad(g\i)X) = (Ad(g\i)'\mu(x))(X).
\endalign$$

\therosteritem5. Let $f\in C^\infty(\g',\Bbb R)$, then we have
$$\align d(\mu^*f)(x)y &= d(f\o\mu)(x)y = df(\mu(x))d\mu(x)y \tag b \\
&= (d\mu(x)y)(df(\mu(x))) = \ome(\rho'(df(\mu(x)))x,y)
\endalign$$ 
by \thetag a, which is smooth in $x$ as a mapping into 
$\bold H_\infty\cong \bold H_\infty^*\subset\bold H_\infty'$
since $\g'$ is finite dimensional. 
 From lemma \nmb!{7.4} we have that
$f\o\mu\in C^\infty_*(\bold H_\infty,\Bbb R)$.

$$\ome(\grad(\mu^*f)(x),y) = d(\mu^*f)(x)y =
\ome(\rho'(df(\mu(x)))x,y)$$
by \thetag b,
so $\grad(\mu^*f)(x) = \rho'(df(\mu(x)))x$. 
The Poisson structure
on $\g'$ is given as follows.
We view the Lie bracket on $\g$ as a linear mapping
$\La^2\g\to\g$; its adjoint $P:\g'\to\La^2\g'$ is then a section
of the bundle $\La^2T\g'\to \g'$, which is called the Poisson
structure on $\g'$. If for $\al\in\g'$ we
view $df(\al)\in L(\g',\Bbb R)$ as an element in $\g$, the
Poisson bracket for $f_i\in C^\infty(\g',\Bbb R)$ is given by
$\{f_1,f_2\}_{\g'}(\al) = (df_1\wedge df_2)(P)|_\al 
= \al([df_1(\al),df_2(\al)])$. Then we may compute as follows.
$$\alignat2 &(\mu^*\{f_1,f_2\}_{\g'})(x) = \{f_1,f_2\}_{\g'}(\mu(x))\\
&\qquad= \mu(x)([df_1(\mu(x)),df_2(\mu(x))])\\
&\qquad= \si([df_1(\mu(x)),df_2(\mu(x))])(x)\\
&\qquad= \{\si(df_1(\mu(x))),\si(df_2(\mu(x)))\}(x)
    &&\text{ by lemma \nmb!{8.2}}\\
&\qquad= \ome(\grad \si(df_2(\mu(x)))(x),\grad \si(df_1(\mu(x)))(x))\\
&\qquad= \ome(\rho'(df_2(\mu(x)))x,\rho'(df_1(\mu(x)))x)\\
&\qquad= \ome(\grad(\mu^*f_2)(x),\grad(\mu^*f_1)(x))
    &&\text{ by \thetag b}\\
&\qquad= \{\mu^*f_1,\mu^*f_2\}_{\bold H_\infty}(x).\qed
\endalignat$$
\enddemo

\heading\totoc \nmb0{9}. The real analytic moment mapping \endheading

\subheading{\nmb.{9.1}} Let again $\rho: G\to U(\bold H)$ be a
unitary representation of a Lie group $G$ on a Hilbert space
$\bold H$. 

\demo{Definition} A vector $x\in \bold H$ is called {it real
analytic} if the mapping $g\mapsto \rho(g)x$, $G\to \bold H$ is
a real analytic mapping, in the real analytic structure of the
Lie group $G$, in the setting explained in section \nmb!{4}. 
\enddemo

Let $\bold H_\a$ denote the vector space of all
real analytic vectors in $\bold H$. 
Then we have a linear
embedding $j:\bold H_\a\to C^\a(G,\bold H)$ into the space of real
analytic mappings, given by 
$x\mapsto (g\mapsto \rho(g)x)$. We equip $C^\a(G,\bold H)$ with
the convenient vector space structure described in \cite{K-M,
5.4, see also 3.13}. 
Then $\bold H_\a$ consists of all equivariant functions in
$C^\a(G,\bold H)$ and is therefore a closed subspace. So it is
a convenient vector space with the induced structure.

The space $\bold H_\a$ is dense in the Hilbert space $\bold H$
by \cite{Wa, 4.4.5.7} and an invariant subspace, so we have a
representation $\rho: G\to L(\bold H_\a,\bold H_\a)$.

\proclaim{\nmb.{9.2}. Theorem} The mapping 
$\hat\rho:G\x\bold H_\a\to \bold H_\a$ is real analytic in the
sense of \rm \cite {K-M}.
\endproclaim
\demo{Proof} By cartesian closedness of the calculus \nmb!{4.6} 
it suffices to show that the canonically associated mapping 
$$\hat\rho\spcheck:G\to C^\a(\bold H_\a,\bold H_\a)$$
is real analytic. It takes values in the closed linear
subspace $L(\bold H_\a,\bold H_\a)$ of all bounded linear
operators. So it suffices to check that the mapping 
$\rho:G\to L(\bold H_\a,\bold H_\a)$ is real analytic. Since
$\bold H_\a$ is a convenient space, by \nmb!{4.5} it
suffices to show that 
$$G @>\rho>> L(\bold H_\a,\bold H_\a) @>{\operatorname{ev}_x}>> \bold H_\a$$
is real analytic for each $x\in \bold H_\a$. Since the
structure on $\bold H_\a$ is induced by the embedding into
$C^\a(G,\bold H)$, we have to check, that 
$$\align &G @>\rho >> L(\bold H_\a,\bold H_\a) @>{\operatorname{ev}_x}>>
\bold H_\a @>j>> C^\a(G,\bold H),\\
&g\mapsto \rho(g)\mapsto \rho(g)x \mapsto (h\mapsto \rho(h)\rho(g)x),
\endalign$$
is real analytic for each $x\in \bold H_\a$. Again by cartesian
closedness \nmb!{4.6} it suffices that the associated mapping
$$\align &G\x G\to \bold H\\
&(g,h)\mapsto \rho(h)\rho(g)x=\rho(hg)x
\endalign$$
is real analytic. And this is the case since $x$ is a real
analytic vector.
\qed\enddemo

\subheading{\nmb.{9.3}} Again we consider now $\bold H_\a$ as a
"weak" symplectic real analytic Fr\'echet manifold, equipped
with the symplectic structure $\ome$, the restriction of the
imaginary part of the hermitian inner product $\langle
\quad,\quad\rangle$ on $\bold H$. Then again $\ome\in \Om^2(\bold
H_\a)$  is a closed 2-form which is non degenerate in the sense that 
$\check\ome:\bold H_\a\to \bold H_\a'=L(\bold H_\a,\Bbb R)$ 
is injective. Let 
$$\bold H_\a^*:= \check\ome(\bold H_\a)=\ome(\bold H_\a,\quad)=
Re\langle \bold H_\a,\quad\rangle\subset \bold H_\a'=L(\bold
H_\a,\Bbb R)$$
again denote the {\it analytic dual} of $\bold H_\a$, equipped
with the topology induced by the isomorphism
with $\bold H_\a$.

\subheading{\nmb.{9.4} Remark} All the results leading to the
smooth moment mapping can now be carried over to the real
analytic setting with {\it no} changes in the proofs. So all
statements from \nmb!{7.5} to \nmb!{8.4} are valid in the real
analytic situation. We summarize this in one more result:

\proclaim{\nmb.{9.5} Theorem} Consider the injective linear
continuous $G$-equivariant mapping $i:\bold H_\a\to \bold
H_\infty$. Then for the smooth moment mapping 
$\mu:\bold H_\infty\to \g'$ from \nmb!{8.4} the composition
$\mu\o i:\bold H_\a\to \bold H_\infty\to\g'$ is real analytic.
It is called the real analytic moment mapping.
\endproclaim
\demo{Proof} It is immediately clear from \nmb!{9.2} and the
formula \nmb!{8.3} for the smooth moment mapping, that $\mu\o i$
is real analytic. 
\qed\enddemo

\enddocument